\documentclass[notitlepage]{article}

 \usepackage[matrix,arrow,curve]{xy}
 \usepackage{graphicx}
 \usepackage[english]{babel}
 \usepackage{amsmath}
 \usepackage{amssymb, indentfirst}
 \usepackage{amsthm}
\begin{document}

\newtheorem{sled}{Corrolary}
\newtheorem{lem}{Lemma}
\newtheorem{zam}{Remark}
\newtheorem{ex}{Example}
\newtheorem{opr}{Definition}
\newtheorem{predl}{Proposition}
\newtheorem{thm}{Theorem}

\title{TWO-TERM TILTING COMPLEXES OVER BRAUER TREE ALGEBRAS}
\author{Alexandra Zvonareva}
\date{}
\maketitle

\begin{abstract}

In this paper all two-term tilting complexes over a Brauer tree
algebra with multiplicity one are described using a classification
of indecomposable two-term partial tilting complexes obtained
earlier in a joint paper with M.~Antipov. The endomorphism rings of
such complexes are computed.

\end{abstract}

\section{Introduction}

Let $A$ be a Brauer tree algebra corresponding to a Brauer tree
$\Gamma$ with multiplicity one. $\text{TrPic}(A)$ is the derived
Picard group of $A$, that is the group of standard autoequivalences
of the derived category of $A$ modulo the natural isomorphism. Let
us consider the derived Picard groupoid, whose objects are the
Brauer tree algebras corresponding to the Brauer trees with $n$
edges, and the morphisms are the standard equivalences between them.
$\text{TrPic}(A)$ is the group of endomorphisms of the object $A$ in
this category. The computation of the derived Picard groupoid seems
to be an easier problem than the computation of $\text{TrPic}(A)$.
The derived Picard group is completely computed only in the case of
algebra with two simple modules \cite{RZ}. In other cases only the
action of different braid groups on $\text{TrPic}(A)$ is known
\cite{RZ}, \cite{IM}, \cite{SI}. On the other hand by the result of
Abe and Hoshino \cite{AH} the derived Picard groupoid corresponding
to the class of Brauer tree algebras with multiplicity of the
exceptional vertex $k$ and a fixed number of simple modules is
generated by one-term and two-term tilting complexes. Thus if we
describe all two-term tilting complexes over $A$, we will obtain the
generating set of the derived Picard groupoid.

The computation of the derived Picard groupoid was the main
motivation while writing this paper. However, the two-term tilting
complexes are also related to $\tau$-tilting theory \cite{AIR} and
to simple-minded systems \cite{Ch}.

This paper is a continuation of the joint paper with M.~Antipov
\cite{AZ}, in which we classified all the indecomposable two-term
partial tilting complexes over Brauer tree algebras with
multiplicity one. In section 3 all two-term tilting complexes are
described in combinatorial terms (theorem 1), in section 4 their
endomorphism rings are computed.

\textbf{Acknowledgement:} I would like to thank Mikhail Antipov for
long and fruitful discussions.

\section{Preliminaries}

Let $K$ be an algebraically closed field, $A$ be a finite
dimensional algebra over $K$. We will denote by $A\text{-}{\rm mod}$
the category of finitely generated left $A$-modules, by $K^b(A)$ --
the bounded homotopy category and by $D^b(A)$ the bounded derived
category of $A\text{-}{\rm mod}.$ The shift functor on the derived
category will be denoted by $[1].$ Let us denote by $A\text{-}{\rm
perf}$ the full subcategory of $D^b(A)$ consisting of perfect
complexes, i.e. of bounded complexes of finitely generated
projective $A$-modules. In the path algebra of a quiver the product
of arrows $\overset{a}{\rightarrow} \overset{b}{\rightarrow}$ will
be denoted by $ab.$

\begin{opr}
A complex $T \in A\text{-}{\rm perf}$ is called tilting if
\begin{enumerate}
    \item $\emph{Hom}_{D^b(A)}(T,T[i])=0, \mbox{ for } i \neq 0$;
    \item T \mbox{ generates }$A\text{-}{\rm perf}  \mbox{ as a triangulated category.}$
\end{enumerate}
\end{opr}

\begin{opr}
A complex $T \in A\text{-}{\rm perf}$ is called partial tilting if
the condition $\emph{1}$ from definition $\emph{1}$ is satisfied.
\end{opr}

\begin{opr}
A tilting complex $T \in A\text{-}{\rm perf}$ is called basic if it
does not contain isomorphic direct summands or equally if
$\emph{End}_{D^b(A)}(T)$ is a basic algebra.
\end{opr}

We will call a (partial) tilting complex a two-term (partial)
tilting complex if it is concentrated in two neighboring degrees.

\begin{opr}

Let $\Gamma$ be a tree with $n$ edges and a distinguished vertex,
which has an assigned multiplicity $k \in \mathbb{N}$ (this vertex
is called exceptional, $k$ is called the multiplicity of the
exceptional vertex). Let us fix a cyclic ordering of the edges
adjacent to each vertex in $\Gamma$ (if $\Gamma$ is embedded into
plane we will assume that the cyclic ordering is clockwise). In this
case $\Gamma$ is called a Brauer tree of type $(n,k)$.

\end{opr}

To a Brauer tree of type $(n,k)$ one can associate an algebra
$A(n,k)$. The algebra $A(n,k)$ is a path algebra of a quiver with
relations. Let us construct a Brauer quiver $Q_{\Gamma}$ using the
Brauer tree $\Gamma$. The vertices of $Q_{\Gamma}$ are the edges of
$\Gamma$, if two edges $i$ and $j$ are incident to the same vertex
in $\Gamma$ and $j$ follows $i$ in the cyclic order of the edges
incident to their common vertex, then there is an arrow from the
vertex $i$ to the vertex $j$ in $Q_{\Gamma}$. $Q_{\Gamma}$ has the
following property: $Q_{\Gamma}$ is the union of oriented cycles
corresponding to the vertices of $\Gamma$, each vertex of
$Q_{\Gamma}$ belongs to exactly two cycles. The cycle corresponding
to the exceptional vertex is called exceptional. The arrows of
$Q_{\Gamma}$ can be divided into two families $\alpha$ and $\beta$
in such a manner that the arrows belonging to intersecting cycles
are in different families.

\begin{opr}
The basic Brauer tree algebra $A(n,k)$, corresponding to a tree
$\Gamma$ of type $(n,k)$ is isomorphic to $Q_{\Gamma}/I$, where the
ideal $I$ is generated by the relations:
\begin{enumerate}
  \item $\alpha\beta=0=\beta\alpha$;
  \item for any vertex
  $x$, not belonging to the exceptional cycle,
  $\alpha^{x_{\alpha}}=\beta^{x_{\beta}}$, where $x_{\alpha},$ resp.
  $x_{\beta}$ is the length of the $\alpha$, resp. $\beta$-cycle,
  containing
  $x$;
  \item for any vertex
  $x$, belonging to the exceptional $\alpha$-cycle (resp. $\beta$-cycle),
  $(\alpha^{x_{\alpha}})^k=\beta^{x_{\beta}}$ (resp.
  $\alpha^{x_{\alpha}}=(\beta^{x_{\beta}})^k$).
\end{enumerate}

An algebra is called a Brauer tree algebra of type $(n,k)$, if it is
Morita equivalent to the algebra $A(n,k)$.

\end{opr}

Note that the ideal $I$ is not admissible. From now on for
convenience all algebras are supposed to be basic.

Rickard showed that two Brauer tree algebras corresponding to the
trees $\Gamma$ and $\Gamma'$ are derived equivalent if and only if
their types $(n,k)$ and $(n',k')$ coincide \cite{Ri2} and it follows
from the results of Gabriel and Riedtmann that this class is closed
under derived equivalence \cite{GR}.

\begin{opr}

Let $B$ be a Brauer tree algebra. $A$-cycle is a maximal ordered set
of nonrepeating arrows of $Q$ such that the product of any two
neighboring arrows is not equal to zero.

\end{opr}

Note that the fact that the algebra is symmetric means that
$A$-cycles are actually cycles.

In \cite{AZ} we classified all indecomposable two-term partial
tilting complexes over a Brauer tree algebra with multiplicity one.
Note that any two-term indecomposable complex is either isomorphic
to a stalk complex of a projective module concentrated in some
degree (such complexes are obviously partial tilting), the minimal
projective presentation of some indecomposable $A$-module.

\textbf{Theorem}\emph{ Let $A$ be a Brauer tree algebra with
multiplicity one. The  minimal projective presentation of an
indecomposable non-projective $A$-module $M$ is a partial tilting
complex if and only if  $M$ is not isomorphic to $P/ \emph{soc}(P)$
for any indecomposable projective module $P.$}

\section{Two-term tilting complexes over Brauer tree algebras with multiplicity one}

To check the conditions from the definitions 2 it will be convenient
to work with partial tilting complexes and not with the
corresponding modules. Let us associate the following diagram on the
quiver of $A$ to a two term complex consisting of projective
modules.

\begin{opr}

Let $T=P_0 \overset{f}{\rightarrow} P_1 \in A-\emph{perf}$ be the
projective presentation of an indecomposable non-projective module
$M,$ $P_0=\bigoplus_{i \in I} Ae_i,$ $P_1=\bigoplus_{i \in J} Ae_i.$
Mark the vertices corresponding to the set $I \cup J$ on the quiver
of $A$. Note that since $M$ does not have repeating composition
factors and since $M$ is not isomorphic to $P/\text{soc}(P)$ for any
indecomposable projective module $P,$ each index can occur only once
and only in one of the sets. Let us also mark the path from $i$ to
$j,$ $i \in I,$ $j \in J$ if $f$ has a nonzero component between the
corresponding summands of $P_0$ and $P_1$ on the quiver of algebra
$A.$ The diagram obtained in such a manner will be called a diagram
of a projective presentation $T=P_0 \overset{f}{\rightarrow} P_1.$
We will call the vertices corresponding to the set $I \cup J$ marked
vertices of the diagram.

\end{opr}
The obtained diagram is a connected path without self-intersections,
it changes its orientation and $A$-cycle in every vertex from the
set $I \cup J.$ This is satisfied since for any index $j \in J$
there exists at most two indices from $I$ such that the
corresponding components of $f$ are nonzero and visa versa: for any
index $i \in I$ there exists at most two indices from $J$ such that
the corresponding components of $f$ are nonzero. It is clear that to
any connected path $\Theta$ without self-intersections, which
changes its orientation every time it changes $A$-cycle and which
consists of more than one vertex, one can associate a two-term
partial tilting complex as follows: let $I$ be the set of indices
corresponding to the sources of $\Theta,$ let $J$ be the set of
indices corresponding to the sinks of $\Theta.$ Then
$P_0=\bigoplus_{i \in I} Ae_i,$ $P_1=\bigoplus_{i \in J} Ae_i,$ $f$
is induced by the morphisms corresponding to the directed subpaths
from the sources to the sinks. (The projective modules corresponding
to the neighboring sink and source belong to the same $A$-cycle and
up to an invertible constant there is a unique morphism between
them, the choice of the coefficient does not play any role, so we
will assume that we always choose the multiplication by the
corresponding path as a morphism.) So there is a one to one
correspondence between minimal projective presentations of
indecomposable non-projective modules non isomorphic to $P/
\emph{soc}(P)$ for any indecomposable projective module $P$ and
connected paths $\Theta$ on the quiver of $A$ without
self-intersections, which change their orientation every time they
change $A$-cycle and which consist of more than one vertex.

For a diagram of a projective presentation and the projective
presentation itself we will often use the same notation.

\begin{opr}

Let $T_i$ and $T_j$ be two diagrams of projective presentations
which meet at more than one vertex so that the intersection has only
one connected component. The restriction of $T_i$ with respect to
$T_j$ is the intersection of the diagram $T_i$ and the union of that
$A$-cycles of the algebra which contain at least one marked vertex
of $T_j$.

\end{opr}

A diagram which consists of more than one vertex and which is
contained in one $A$-cycle will be called a string. As in \cite{AZ}
such a diagram will be denoted by $(k,...,l)$, where $k$ is a sink,
and $l$ is a source of the string.

\begin{zam}
Thus the restriction of $T_i$ with respect to $T_j$ is some
subdiagram of $T_i$ which contains the intersection of $T_i$ and
$T_j$ and the completion of the intersection to the restriction can
be defined independently at the ends of the intersection. If the
intersection of $T_i$ and $T_j$ does not contain a marked vertex of
$T_i$ the restriction is a substring of $T_i$ which contains the
intersection; if the intersection ends at an unmarked vertex of
$T_i$ let us complete it to the smallest subdiagram (to the nearest
marked vertex); if the intersection ends at a marked vertex of $T_i$
which is not a marked vertex of $T_j$ we will not complete it, if
the intersection ends at a marked vertex of $T_i$ which is a marked
vertex of $T_j$, let us complete it to the next marked vertex of
$T_i$, if it exists.
\end{zam}

\begin{zam}
The restriction is defined in such a manner that there are no
nonzero morphisms between the projective summands of the components
of $T_i$ and $T_j$ which correspond to the vertices not contained in
the restrictions.
\end{zam}

To classify all basic two-term tilting complexes it is necessary and
sufficient to classify $n$-tuples of pairwise orthogonal
nonisomorphic indecomposable two-term partial tilting complexes
\cite{AH}.

\begin{thm}
Let $T_i$ and $T_j$ be indecomposable partial tilting complexes. The
complex $T_i \oplus T_j$ is partial tilting iff one of the following
conditions holds:

1) In the case when $T_i$ and $T_j$ are projective presentations of
some modules, such that $T_i$ and $T_j$ are indecomposable two-term
partial tilting complexes.

a) The diagrams of $T_i$ and $T_j$ do not have vertices belonging to
the same $A$-cycle.

b) The diagrams of $T_i$ and $T_j$ have vertices belonging to the
same $A$-cycle $\Upsilon$, but they do not intersect, and $\Upsilon$
does not contain a source of degree one of $T_j$ (resp. $T_i$) and a
sink of degree one of $T_i$ (resp. $T_j$) such that these are the
only vertices of $T_i$ and $T_j$ belonging to $\Upsilon$.

c) The diagrams of $T_i$ and $T_j$ meet at one vertex $k$,and $k$ is
neither a marked vertex of $T_i$ nor a marked vertex of $T_j,$ or
$k$ is a degree one sink of both $T_i$ and $T_j,$ or $k$ is a degree
one source of both $T_i$ and $T_j.$

d) The intersection of $T_i$ and $T_j$ consists of one connected
component, which contains more than one vertex. The diagrams of
$T_i$ and $T_j$ intersect in such a way that one of the end points
of the intersection is a sink and another is a source, the
restriction of $T_i$ with respect to $T_j$ belongs to the
restriction of $T_j$ with respect to $T_i$ or visa versa.

e) The intersection of $T_i$ and $T_j$ consists of one connected
component, which contains more than one vertex. The diagrams of
$T_i$ and $T_j$ intersect in such a way that both of the end
vertices of the intersection are sinks and neither the restriction
of $T_i$ with respect to $T_j$ belongs to the restriction of $T_j$
with respect to $T_i$ nor the restriction of $T_j$ with respect to
$T_i$ belongs to the restriction of $T_i$ with respect to $T_j$; or
the diagrams of $T_i$ and $T_j$ have a coinciding vertex of degree
one.

f) The intersection of $T_i$ and $T_j$ consists of one connected
component, which contains more than one vertex. The diagrams of
$T_i$ and $T_j$ intersect in such a way that both of the end
vertices of the intersection are sources and neither the restriction
of $T_i$ with respect to $T_j$ belongs to the restriction of $T_j$
with respect to $T_i$ nor the restriction of $T_j$ with respect to
$T_i$ belongs to the restriction of $T_i$ with respect to $T_j$; or
the diagrams of $T_i$ and $T_j$ have a coinciding vertex of degree
one.

2) In the case when $T_i$ is the projective presentation of some
module, such that $T_i$ is an indecomposable two-term partial
tilting complex and $T_j$ is an indecomposable stalk complex of a
projective module $P$.

a) $P$ is concentrated in 0, the vertex corresponding to $P$
coincides with a source of degree one of $T_i$, or the vertex
corresponding to $P$ does not belong to $T_i$ and there is no
$A$-cycle which contains the vertex corresponding to $P$ and a sink
of degree one of $T_i$ such that this is the only vertex of $T_i$
belonging to this $A$-cycle.

b) $P$ is concentrated in 1, the vertex corresponding to $P$
coincides with a sink of degree one of $T_i$, or the vertex
corresponding to $P$ does not belong to $T_i$ and there is no
$A$-cycle which contains the vertex corresponding to $P$ and a
source of degree one of $T_i$ such that this is the only vertex of
$T_i$ belonging to this $A$-cycle.

3) In the case when $T_i$ and $T_j$ are two indecomposable stalk
complexes of projective modules and the vertices corresponding to
$T_i$ and $T_j$ do not belong to the same $A$-cycle, or the vertices
corresponding to $T_i$ and $T_j$ belong to the same $A$-cycle, $T_i$
and $T_j$ are concentrated in the same degree.

\end{thm}

The rest of section 3 is devoted to the proof of this theorem.

As mentioned before, to classify all basic two-term tilting
complexes it is necessary and sufficient to classify $n$-tuples of
nonisomorphic indecomposable two-term partial tilting complexes
$\{T_1,...,T_n\}$ such that
$\emph{Hom}_{D^b(A)}(T_i,T_j[1])=0=\emph{Hom}_{D^b(A)}(T_i,T_j[-1])$
\cite{AH}. By the following remark by Happel \cite{Ha} it is
sufficient to check only one of these conditions, for example,
$\emph{Hom}_{D^b(A)}(T_i,T_j[-1]).$

\begin{zam}

Let $A$ be a finite dimensional algebra over a field $K$, let
$\emph{proj}-A$, $\emph{inj}-A$ be the categories of finitely
generated projective and injective modules, $K^b(\emph{proj}-A)$,
$K^b(\emph{inj}-A)$ the corresponding bounded homotopy categories,
$D$ the duality with respect to $K.$ Then the Nakayama functor $\nu$
induces an equivalence of triangulated categories
$K^b(\emph{proj}-A) \rightarrow K^b(\emph{inj}-A)$ and there is a
natural isomorphism $D\emph{Hom}(P,-)\rightarrow \emph{Hom}(-,\nu
P),$ $P \in K^b(\emph{proj}-A)$.

\end{zam}

First, let us consider the case of two projective presentations of
some modules.

\begin{lem}
Let $T_i$ and $T_j$ not have vertices belonging to the same
$A$-cycle then $\emph{Hom}_{D^b(A)}(T_i,T_j[-1])=0.$

\end{lem}

\textbf{Proof.} This is obvious since there are no nonzero morphisms
between the summand of the component of $T_i$ and
$T_j$.~\hfill\(\Box\)

\begin{lem}
Let the diagrams of $T_i$ and $T_j$ have vertices belonging to the
same $A$-cycle $\Upsilon$, but not intersect, then the condition
$\emph{Hom}_{D^b(A)}(T_i,T_j[-1])=0=\emph{Hom}_{D^b(A)}(T_j,T_i[-1])$
holds in all the cases but the case when $\Upsilon$ contains a
source of degree one of $T_j$ (resp. $T_i$) and a sink of degree one
of $T_i$ (resp. $T_j$) such that these are the only vertices of
$T_i$ and $T_j$ belonging to $\Upsilon$.
\end{lem}

\textbf{Proof.} Note that if the diagrams of $T_i$ and $T_j$ have
vertices belonging to the same $A$-cycle $\Upsilon$, but do not
intersect, then since $\Gamma$ is a tree there is no other $A$-cycle
they both meet. This situation can occur in the following cases:
$\Upsilon$ contains a substring of $T_i$ and a substring of $T_j$, a
substring of $T_i$ and a vertex of $T_j$, a vertex of $T_i$ and a
vertex of $T_j$, namely:

1) $\Upsilon$ contains a substring $(i_1,...,i_2)$ of $T_i$ and a
substring $(j_1,...,j_2)$ of $T_j$ but they do not intersect. Let us
show that there is no chain morphism between $T_i$ and $T_j[-1]$. It
is clear that there is a nonzero morphism only between the
projective summands of the components of $T_i$ and $T_j$ with
corresponding vertices belonging to $\Upsilon$ (i.e. between
$Ae_{i_1}$, $Ae_{i_2}$, $Ae_{j_1}$, $Ae_{j_2}$), since all other
projective summands of the components of $T_i$ and $T_j$ belong to
different $A$-cycles. It is clear that
$\emph{Hom}_{D^b(A)}(T_i,T_j[-1])=0$ iff
$\emph{Hom}_{D^b(A)}(Ae_{i_2} \rightarrow Ae_{i_1},Ae_{j_2}
\rightarrow Ae_{j_1}[-1])=0,$ which is true since the strings
$(i_1,...,i_2)$ and $(j_1,...,j_2)$ do not intersect.

2) $\Upsilon$ contains a substring $(i_1,...,i_2)$ of $T_i$ and a
vertex of $T_j$ which is not marked.
$\emph{Hom}_{D^b(A)}(T_i,T_j[-1])=0$, since the projective summands
of the components of $T_i$ and $T_j$ belong to different $A$-cycles.
$\emph{Hom}_{D^b(A)}(T_j,T_i[-1])=0$ for the same reasons.

3) $\Upsilon$ contains a substring $(i_1,...,i_2)$ of $T_i$ and a
vertex $(j_1)$ of $T_j$ which is a sink of degree one. Then there is
no nonzero morphism from $Ae_{i_1}$ to the component of $T_j[-1],$
concentrated in $1,$ i.e. $\emph{Hom}_{D^b(A)}(T_i,T_j[-1])=0$. It
is clear that $\emph{Hom}_{D^b(A)}(T_j,T_i[-1])=0$ iff
$\emph{Hom}_{D^b(A)}(0 {\rightarrow} Ae_{j_1},Ae_{i_2} \rightarrow
Ae_{i_1}[-1])=0,$ which is true since the string $(i_1,...,i_2)$
does not contain $(j_1)$.

4) $\Upsilon$ contains a substring $(i_1,...,i_2)$ of $T_i$ and a
vertex $(j_1)$ of $T_j$ which is a source of degree one. This case
is similar to the previous one.

5) $\Upsilon$ contains a vertex $(i_1)$ of $T_i$ and a vertex
$(j_1)$ of $T_j$, both vertices are sources of degree one and
$\Upsilon$ does not contain any other vertices of $T_i$ and $T_j$.
The modules $Ae_{i_1}$ and $Ae_{j_1}$ are the only projective
summands of the components of $T_i$ and $T_j$ with a nonzero
hom-space, so after we apply a shift functor there will be no
nonzero morphisms.

6) $\Upsilon$ contains a vertex $(i_1)$ of $T_i$ and a vertex
$(j_1)$ of $T_j$, both vertices are sinks of degree one and
$\Upsilon$ does not contain any other vertices of $T_i$ and $T_j$.
This case is similar to the previous one.

7) $\Upsilon$ contains a degree one sink of $T_i$ and an unmarked
vertex of $T_j$ and $\Upsilon$ does not contain any other vertices
of $T_i$ and $T_j$. $\emph{Hom}_{D^b(A)}(T_i,T_j[-1])=0$, since the
projective summands of the components of $T_i$ and $T_j$ belong to
different $A$-cycles. $\emph{Hom}_{D^b(A)}(T_j,T_i[-1])=0$ for the
same reasons.

8) $\Upsilon$ contains a degree one source of $T_i$ and an unmarked
vertex of $T_j$ and $\Upsilon$ does not contain any other vertices
of $T_i$ and $T_j$. This case is similar to the previous one.

9) $\Upsilon$ contains an unmarked vertex of $T_i$ and an unmarked
vertex of $T_j$ and $\Upsilon$ does not contain any other vertices
of $T_i$ and $T_j$. This case is similar to the previous one.

10) $\Upsilon$ contains a degree one sink $(i_1)$ of $T_i$ and a
degree one source $(j_1)$ of $T_j$. In this case $Ae_{i_1}$ and
$Ae_{j_1}$ are projective summands of components of $T_i$ and $T_j,$
which are concentrated in different degrees, after we apply the
shift functor they will be concentrated in the same degree and so
there will be a nonzero morphism between them, which induces a chain
map. Hence $\emph{Hom}_{D^b(A)}(T_i,T_j[-1]) \neq 0$.~\hfill\(\Box\)

\begin{lem}

Let diagrams of $T_i$ and $T_j$ have an intersection containing more
than one connected component, then at least one of the spaces
$\emph{Hom}_{D^b(A)}(T_i,T_j[-1]),$
$\emph{Hom}_{D^b(A)}(T_j,T_i[-1])$ is nonzero.

\end{lem}

\textbf{Proof.} Let diagrams of $T_i$ and $T_j$ have an intersection
containing at least 3 connected components, then one of these
components is an isolated vertex $k$ (to get from one $A$-cycle to
another one should pass a unique vertex and since the diagrams have
an intersection containing at least 3 connected components near this
vertex they have an opposite orientation). Hence for one of the
diagrams this vertex is a sink (say for $T_i$) and for the other
this vertex is a source (for $T_j$). There is a nonzero morphism
from $Ae_k$ to itself whose image is the socle of $Ae_k$. It induces
a nonzero morphism from $\emph{Hom}_{D^b(A)}(T_i,T_j[-1]).$

Let diagrams of $T_i$ and $T_j$ have an intersection containing two
connected components, it is clear that both components belong to the
same $A$-cycle. Let $(i_1,...,i_2)$, $(j_1,...,j_2)$ be substrings
of $T_i$ and $T_j$ corresponding to the restrictions of $T_i$ and
$T_j$ to this $A$-cycle. The vertices have the following  order on
the $A$-cycle: $i_1,i_2,j_1,j_2$, where $i_1$ can coincide with
$j_2$ and $j_1$ can coincide with $i_2.$ There is a nonzero morphism
from $Ae_{i_1}$ to $Ae_{j_2}$ (in the case when $i_1$ coincides with
$j_2$ the morphism which has the socle as its image), it induces a
chain map.~\hfill\(\Box\)

\begin{lem}

The diagrams of $T_i$ and $T_j$ have an intersection containing one
vertex $k$, then
$\emph{Hom}_{D^b(A)}(T_i,T_j[-1])=0=\emph{Hom}_{D^b(A)}(T_j,T_i[-1])$
holds in the following cases: $k$ is neither a marked vertex of
$T_i$ nor a marked vertex of $T_j,$ or $k$ is a sink of degree one
of both $T_i$ and $T_j,$ or $k$ is a source of degree one of both
$T_i$ and $T_j.$

\end{lem}

\textbf{Proof.} Note that there is no $A$-cycle which intersects
with $T_i$ and $T_j$ but does not contain $k$ (since $\Gamma$ is a
tree). The vertex $k$ can be a degree one sink, a degree one source
or an unmarked vertex of both $T_i$ and $T_j$, so we have to
consider 6 cases.

1) Let $k$ be an unmarked vertex of both $T_i$ and $T_j$ then marked
vertices of $T_i$ and marked vertices of $T_j$ belong to different
$A$-cycles there are no nonzero morphisms between the corresponding
projective modules.

2) Let $k$ be a degree one sink of both $T_i$ and $T_j$ then
$\emph{Hom}_{D^b(A)}(T_i,T_j[-1])=0$ iff
$\emph{Hom}_{D^b(A)}(Ae_{i_2}\rightarrow Ae_k,(Ae_{j_2}\rightarrow
Ae_k)[-1])=0$ (where $i_2$, $j_2$ are the marked vertices of $T_i$
and $T_j$ which are next to $k$), which is true since the only
nonzero morphism from $\emph{Hom}_A(Ae_k,Ae_{j_2})$ does not induce
a chain map. The space $\emph{Hom}_{D^b(A)}(T_j,T_i[-1])$ is zero
for the same reasons.

3) The case when $k$ is a degree one source of both $T_i$ and $T_j$
is similar to the previous one.

Let us check that in all other cases at least one of the spaces
$\emph{Hom}_{D^b(A)}(T_i,T_j[-1])$,
$\emph{Hom}_{D^b(A)}(T_j,T_i[-1])$ is not zero.

4) The case when $k$ is a source of $T_i$ and a sink of $T_j$: the
morphism $Ae_k\rightarrow Ae_k$ which has a socle of $Ae_k$ as its
image annihilates any noninvertible morphism between indecomposable
projective modules, thus it induces a nonzero chain map from
$\emph{Hom}_{D^b(A)}(T_j,T_i[-1])$.

5) The case when $k$ is a source of $T_i$ but is not a marked vertex
of $T_j$. Let $(j_1,...,j_2)$ be a substring of $T_j$ containing
$k$, $(i_1,...,k)$ a substring of $T_i$, the corresponding vertices
are ordered in the following way $j_1, j_2, k$ on the $A$-cycle. The
morphism $Ae_{j_1}\rightarrow Ae_k$ induces a nonzero chain map from
$\emph{Hom}_{D^b(A)}(T_j,T_i[-1])$.

6) The case when $k$ is a sink of $T_i$ but is not a marked vertex
of $T_j$ is similar to the previous one.~\hfill\(\Box\)

If the intersection of $T_i$ and $T_j$ consists of more than one
vertex then the following types of intersections can occur: one of
the end points of the intersection is a sink and another is a
source, both of the end points of the intersection are sinks, both
of the end points of the intersection are sources.

\begin{lem}

Let the intersection of $T_i$ and $T_j$ consist of one connected
component, which contains more than one vertex. The diagrams of
$T_i$ and $T_j$ intersect in such a way that one of the end points
of the intersection is a sink and another is a source, then
$\emph{Hom}_{D^b(A)}(T_i,T_j[-1])=0=\emph{Hom}_{D^b(A)}(T_j,T_i[-1])$
iff the restriction of $T_i$ with respect to $T_j$ belongs to the
restriction of $T_j$ with respect to $T_i$ or visa versa.

\end{lem}

\textbf{Proof.} Let us first consider the cases when the end points
of the intersection coincide with the marked vertices of only one of
the diagrams (the cases 1 and 2), after that let us consider the
cases when the end points of the intersection coincide with the
marked vertices of both diagrams (the cases 3, 4, 5).

1) Assume that the restriction of $T_i$ with respect to $T_j$
belongs to the restriction of $T_j$ with respect to $T_i$ and the
end points of the intersection of $T_i$ and $T_j$ do not coincide
with the marked vertices of $T_j$. The corresponding subdiagrams are
arranged on the quiver of $A$ as follows:

$$
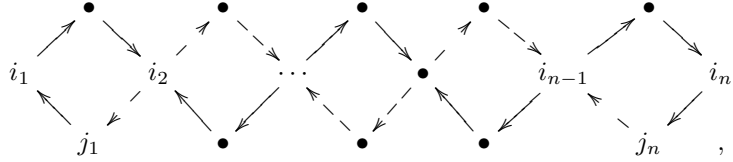
\xymatrix @=1pc{
&\bullet \ar[dr] &&\bullet\ar@{-->}[dr]&&\bullet\ar[dr]&&\bullet\ar@{-->}[dr]&&\bullet\ar[dr]& \\
i_1\ar[ur]&&i_2\ar@{-->}[ur]\ar@{-->}[dl]&&\cdots\ar[ur]\ar[dl]&&\bullet\ar@{-->}[ur]\ar@{-->}[dl]&&i_{n-1}\ar[dl]\ar[ur]&&i_n\ar[dl] \\
&j_1\ar[ul]&&\bullet\ar[ul]&&\bullet\ar@{-->}[ul]&&\bullet\ar[ul]&&j_n\ar@{-->}[ul]&, \\
}
$$
where $i_1, i_2,...,i_n$ are the marked vertices of $T_i$ (which
belong to the restriction), the vertices with odd indices are
sources, the vertices with even indices are sinks; $j_1,
i_2,...,i_n-1,j_n$ are the marked vertices of $T_j$ (which belong to
the restriction), the vertices with odd indices are sources, the
vertices with even indices are sinks ($i_1 \neq j_1, i_n \neq j_n$).
The differentials of the restrictions of $T_i$ and $T_j$ are induced
by the morphisms from $Ae_{i_k}$ to $Ae_{i_{k-1}}, Ae_{i_{k+1}}$ for
odd $k=3,...,n-3$, and the morphisms from $Ae_{i_1}$ to
$Ae_{i_{2}},$ from $Ae_{j_1}$ to $Ae_{i_{2}},$  from $Ae_{i_{n-1}}$
to $Ae_{i_{n-2}}$ and $Ae_{i_{n}},$ from $Ae_{i_{n-1}}$ to
$Ae_{i_{n-2}}$ and $Ae_{j_{n}}$, the corresponding morphisms are
given by the multiplication by the unique nonzero paths on the
quiver of $A$. As it was already mentioned, there are no nonzero
morphisms between the projective summands of the components of $T_i$
and $T_j$ which do not belong to the restrictions.

Assume that there is a nonzero chain map $g$ from the restriction of
$T_j$ to the restriction of $T_i$ shifted by $-1$, it is nonzero on
some projective summand of the component of $T_j$ concentrated in
$1$. Let $g$ restricted to $Ae_{i_k}$ be nonzero, where $k$ is even.
If $g$ has a nonzero component from $Ae_{i_k}$ to $Ae_{i_{k-1}},$
the composition of $g$ and the differential should be equal to 0,
but the composition $Ae_{i_k} \rightarrow Ae_{i_{k-1}} \rightarrow
Ae_{i_k}$ is not equal to zero, hence $g$ has a nonzero component
from $Ae_{i_k}$ to $Ae_{i_{k+1}},$ but then $Ae_{i_{k+1}}
\rightarrow Ae_{i_{k}} \rightarrow Ae_{i_{k+1}}$ is not equal to
zero, hence $g$ has a nonzero component from $Ae_{i_{k+2}}$ to
$Ae_{i_{k+1}},$ therefore $g$ has a nonzero component from
$Ae_{j_n}$ to $Ae_{i_{n-1}},$ but $Ae_{j_{n}} \rightarrow
Ae_{i_{n-1}} \rightarrow Ae_{i_{n}}$ is not equal to zero. Hence if
$g$ has a nonzero component from $Ae_{i_k}$ to $Ae_{i_{k-1}},$ then
$g$ can not be a chain map; if $g$ has a nonzero component from
$Ae_{i_k}$ to $Ae_{i_{k+1}},$ then $g$ can not be a chain map for a
similar reason.

Let us consider $\emph{Hom}_{D^b(A)}(T_i,T_j[-1])$. Assume that
there is a nonzero chain map $g$ from the restriction of $T_i$ to
the restriction of $T_j$ shifted by $-1$, it is nonzero on some
projective summand of the component of $T_i$ concentrated in $1$.
Let $g$ restricted to $Ae_{i_k}$ be nonzero, where $k$ is even. If
$g$ has a nonzero component from $Ae_{i_k}$ to $Ae_{i_{k-1}},$ the
composition of $g$ and the differential should be equal to 0, but
the composition $Ae_{i_{k-1}} \rightarrow Ae_{i_{k}} \rightarrow
Ae_{i_{k-1}}$ is not equal to zero, hence $g$ has a nonzero
component from $Ae_{i_{k-2}}$ to $Ae_{i_{k-1}},$ but then
$Ae_{i_{k-2}} \rightarrow Ae_{i_{k-1}} \rightarrow Ae_{i_{k-2}}$ is
not equal to zero, hence $g$ has a nonzero component from
$Ae_{i_{k-2}}$ to $Ae_{i_{k-3}},$ in such a manner we get that $g$
has a nonzero component from $Ae_{i_2}$ to $Ae_{j_{1}},$ but then
$Ae_{i_{1}} \rightarrow Ae_{i_{2}} \rightarrow Ae_{j_{1}}$ is not
equal to zero. Hence if $g$ has a nonzero component from $Ae_{i_k}$
to $Ae_{i_{k-1}},$ then $g$ can not be a chain map; if $g$ has a
nonzero component from $Ae_{i_k}$ to $Ae_{i_{k+1}},$ then $g$ can
not be a chain map for a similar reason. If $g$ has a nonzero
component from $Ae_{i_2}$ to $Ae_{j_{1}},$ then it is
straightforward that $g$ can not be a chain map.

2) Consider the case when neither the restriction of $T_i$ with
respect to $T_j$ belongs to the restriction of $T_j$ with respect to
$T_i$ nor the restriction of $T_j$ with respect to $T_i$ belongs to
the restriction of $T_i$ with respect to $T_j$, and $i_1 \neq j_1,
i_n \neq j_n$. This is satisfied if the degree one source of the
restriction of $T_i$ does not belong to the intersection and the
degree one sink of $T_i$ belongs to the intersection (or the same
holds for $T_j$, but this case is analogous). The corresponding
subdiagrams are arranged on the quiver of $A$ as follows:

$$
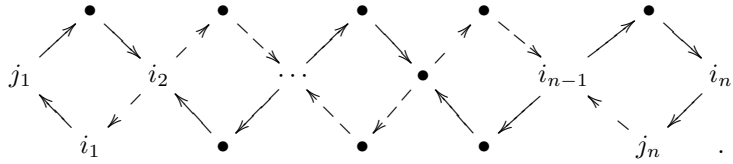
\xymatrix @=1pc{
&\bullet \ar[dr] &&\bullet\ar@{-->}[dr]&&\bullet\ar[dr]&&\bullet\ar@{-->}[dr]&&\bullet\ar[dr]& \\
j_1\ar[ur]&&i_2\ar@{-->}[ur]\ar@{-->}[dl]&&\cdots\ar[ur]\ar[dl]&&\bullet\ar@{-->}[ur]\ar@{-->}[dl]&&i_{n-1}\ar[dl]\ar[ur]&&i_n\ar[dl] \\
&i_1\ar[ul]&&\bullet\ar[ul]&&\bullet\ar@{-->}[ul]&&\bullet\ar[ul]&&j_n\ar@{-->}[ul]&. \\
}
$$

Let us construct a nonzero morphism from
$\emph{Hom}_{D^b(A)}(T_i,T_j[-1])$. Let us describe how this
morphism acts on the components: choose an arbitrary nonzero
morphism $Ae_{i_{2}} \rightarrow Ae_{j_{1}}$, note that the
composition $Ae_{i_{1}} \rightarrow Ae_{i_{2}} \rightarrow
Ae_{j_{1}}$ is equal to zero; choose such a morphism $Ae_{i_{2}}
\rightarrow Ae_{i_{3}}$ that the composition $Ae_{i_{2}} \rightarrow
Ae_{j_{1}} \rightarrow Ae_{i_{2}}$ is equal to $Ae_{i_{2}}
\rightarrow Ae_{i_{3}} \rightarrow Ae_{i_{2}}$ with the opposite
sign; choose such a morphism $Ae_{i_{4}} \rightarrow Ae_{i_{3}}$
that the composition $Ae_{i_{3}} \rightarrow Ae_{i_{2}} \rightarrow
Ae_{i_{3}}$ is equal to the composition $Ae_{i_{3}} \rightarrow
Ae_{i_{4}} \rightarrow Ae_{i_{3}}$ with the opposite sign; continue
to construct the morphism in this way until we reach $Ae_{i_{n}}
\rightarrow Ae_{i_{n-1}}$, which was chosen such that $Ae_{i_{n-1}}
\rightarrow Ae_{i_{n-2}} \rightarrow Ae_{i_{n-1}}$ is equal to the
composition $Ae_{i_{n-1}} \rightarrow Ae_{i_{n}} \rightarrow
Ae_{i_{n-1}}$ with the opposite sign, we can see that the
composition $Ae_{i_{n}} \rightarrow Ae_{i_{n-1}} \rightarrow
Ae_{j_{n}}$ is equal to zero. Thus we have constructed a nonzero
chain map from $\emph{Hom}_{D^b(A)}(T_i,T_j[-1])$.

3) Let now $i_1=j_1$, $i_n \neq j_n,$ assume that $i_n$ belongs to
the restriction of $T_j$. The following options are possible: a) the
vertex $i_1$ is a degree one source of both $T_i$ and $T_j$; b) the
vertex $i_1$ is a degree one source of $T_i$ but not $T_j$ c) the
vertex $i_1$ is a degree one source of  $T_j$ but not $T_i$.

a) The vertex $i_1$ is a degree one source of both $T_i$ and $T_j.$
As before we can see that if there is a nonzero chain map $g$ from
the restriction of $T_j$ to the restriction of $T_i$ shifted by
$-1$, then it should have a nonzero component from $Ae_{i_2}$ to
$Ae_{j_1}=Ae_{i_{1}},$ but $Ae_{i_{1}} \rightarrow Ae_{i_{2}}
\rightarrow Ae_{j_{1}}=Ae_{i_{1}}$ is not equal to zero, hence $g$
can not be a chain map, thus $\emph{Hom}_{D^b(A)}(T_j,T_i[-1])=0.$
It is clear that $\emph{Hom}_{D^b(A)}(T_i,T_j[-1])=0$ for similar
reasons.

b) The vertex $i_1$ is a degree one source of $T_i$ but not $T_j$.
Let $j_0$ be the marked vertex of $T_j$ next to $j_1$. The
restrictions of $T_i$ and $T_j$ are arranged on the quiver of $A$ as
follows:

$$
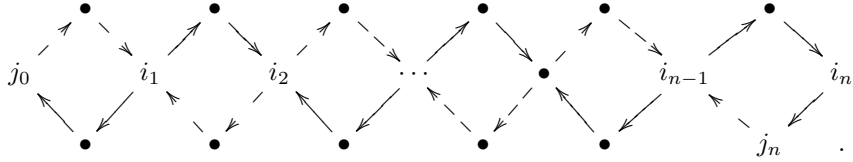
\xymatrix @=1pc{
&\bullet \ar@{-->}[dr] &&\bullet \ar[dr] &&\bullet\ar@{-->}[dr]&&\bullet\ar[dr]&&\bullet\ar@{-->}[dr]&&\bullet\ar[dr]& \\
j_0 \ar@{-->}[ur] &&i_1\ar[ur] \ar[dl]&&i_2\ar@{-->}[ur]\ar@{-->}[dl]&&\cdots\ar[ur]\ar[dl]&&\bullet\ar@{-->}[ur]\ar@{-->}[dl]&&i_{n-1}\ar[dl]\ar[ur]&&i_n\ar[dl] \\
&\bullet \ar[ul] &&\bullet\ar@{-->}[ul]&&\bullet\ar[ul]&&\bullet\ar@{-->}[ul]&&\bullet\ar[ul]&&j_n\ar@{-->}[ul]&. \\
}
$$

As before we can see that if there is a nonzero chain map $g$ from
the restriction of $T_i$ to the restriction of $T_j$ shifted by
$-1$, then it should have a nonzero component from $Ae_{i_2}$ to
$Ae_{i_1},$ but $Ae_{i_{1}} \rightarrow Ae_{i_{2}} \rightarrow
Ae_{i_{1}}$ is not equal to zero, hence $g$ can not be a chain map,
thus $\emph{Hom}_{D^b(A)}(T_i,T_j[-1])=0.$ Let us consider
$\emph{Hom}_{D^b(A)}(T_j,T_i[-1])=0$. If there is a nonzero chain
map $g$ from the restriction of $T_j$ to the restriction of $T_i$
shifted by $-1$, then it should have a nonzero component from
$Ae_{j_n}$ to $Ae_{i_{n-1}},$ but $Ae_{j_{n}} \rightarrow
Ae_{i_{n-1}} \rightarrow Ae_{i_{n}}$ is not equal to zero, hence $g$
can not be a chain map.

c) The vertex $i_1$ is a degree one source of  $T_j$ but not $T_i$.
Let $i_0$ be the marked vertex of $T_i$ next to $i_1$. The
restrictions of $T_i$ and $T_j$ are arranged on the quiver of $A$ as
follows:

$$
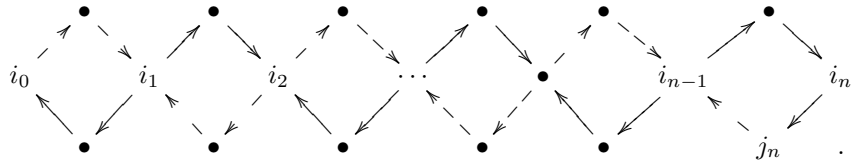
\xymatrix @=1pc{
&\bullet \ar@{-->}[dr] &&\bullet \ar[dr] &&\bullet\ar@{-->}[dr]&&\bullet\ar[dr]&&\bullet\ar@{-->}[dr]&&\bullet\ar[dr]& \\
i_0 \ar@{-->}[ur] &&i_1\ar[ur] \ar[dl]&&i_2\ar@{-->}[ur]\ar@{-->}[dl]&&\cdots\ar[ur]\ar[dl]&&\bullet\ar@{-->}[ur]\ar@{-->}[dl]&&i_{n-1}\ar[dl]\ar[ur]&&i_n\ar[dl] \\
&\bullet \ar[ul] &&\bullet\ar@{-->}[ul]&&\bullet\ar[ul]&&\bullet\ar@{-->}[ul]&&\bullet\ar[ul]&&j_n\ar@{-->}[ul]&. \\
}
$$

Let us construct a nonzero morphism from
$\emph{Hom}_{D^b(A)}(T_i,T_j[-1])$. This morphism is constructed the
same way as in the case when neither the restriction of $T_i$ with
respect to $T_j$ belongs to the restriction of $T_j$ with respect to
$T_i$ nor the restriction of $T_j$ with respect to $T_i$ belongs to
the restriction of $T_i$ with respect to $T_j$, and $i_1 \neq j_1,
i_n \neq j_n$. Let us start from the morphism $Ae_{i_{0}}
\rightarrow Ae_{i_{1}}$, and let us construct the morphism as
described earlier until we reach $Ae_{i_{n}} \rightarrow
Ae_{i_{n-1}}$, which was chosen such that the composition
$Ae_{i_{n-1}} \rightarrow Ae_{i_{n-2}} \rightarrow Ae_{i_{n-1}}$ is
equal to the composition $Ae_{i_{n-1}} \rightarrow Ae_{i_{n}}
\rightarrow Ae_{i_{n-1}}$ with the opposite sign, we can see that
the composition $Ae_{i_{n}} \rightarrow Ae_{i_{n-1}} \rightarrow
Ae_{j_{n}}$ is equal to zero. Thus we have constructed a nonzero
chain map from $\emph{Hom}_{D^b(A)}(T_i,T_j[-1])$.

4) Let $i_n=j_n$, $i_1 \neq j_1,$ assume that $i_1$ belongs to the
restriction of $T_j$. The following options are possible: a) the
vertex $i_n$ is a degree one sink of both $T_i$ and $T_j$; b) the
vertex $i_n$ is a degree one sink of $T_i$ but not $T_j$ c) the
vertex $i_n$ is a degree one sink of  $T_j$ but not $T_i$.

a) The vertex $i_n$ is a degree one sink of both $T_i$ and $T_j.$ A
nonzero chain map from the restriction of $T_i$ to the restriction
of $T_j$ shifted by $-1$, should have a nonzero component from
$Ae_{i_n}$ to $Ae_{i_{n-1}}$, but $Ae_{i_{n}} \rightarrow
Ae_{i_{n-1}} \rightarrow Ae_{j_{n}}=Ae_{i_{n}}$ is not equal to
zero, hence $\emph{Hom}_{D^b(A)}(T_i,T_j[-1])=0.$ For similar
reasons as before $\emph{Hom}_{D^b(A)}(T_j,T_i[-1])=0.$

b) The vertex $i_n$ is a degree one sink of $T_i$ but not $T_j$. Let
$j_{n+1}$ be the marked vertex of $T_j$ next to $i_n$. The
restrictions of $T_i$ and $T_j$ are arranged on the quiver of $A$ as
follows:

$$
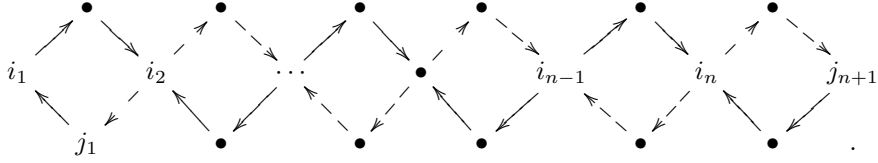
\xymatrix @=1pc{
&\bullet \ar[dr] &&\bullet\ar@{-->}[dr]&&\bullet\ar[dr]&&\bullet\ar@{-->}[dr]&&\bullet\ar[dr]& &\bullet\ar@{-->}[dr]& \\
i_1\ar[ur]&&i_2\ar@{-->}[ur]\ar@{-->}[dl]&&\cdots\ar[ur]\ar[dl]&&\bullet\ar@{-->}[ur]\ar@{-->}[dl]&&i_{n-1}\ar[dl]\ar[ur]&&i_n\ar@{-->}[dl] \ar@{-->}[ur] && j_{n+1}\ar[dl] \\
&j_1\ar[ul]&&\bullet\ar[ul]&&\bullet\ar@{-->}[ul]&&\bullet\ar[ul]&&\bullet\ar@{-->}[ul]& &\bullet\ar[ul]&. \\
}
$$

A nonzero chain map from the restriction of $T_i$ to the restriction
of $T_j$ shifted by $-1$, should have a nonzero component from
$Ae_{i_2}$ to $Ae_{j_1},$ but $Ae_{i_{1}} \rightarrow Ae_{i_{2}}
\rightarrow Ae_{j_{1}}$ is not equal to zero, hence
$\emph{Hom}_{D^b(A)}(T_i,T_j[-1])=0.$ If there is a nonzero chain
map $g$ from the restriction of $T_j$ to the restriction of $T_i$
shifted by $-1$, it should have a nonzero component from $Ae_{i_n}$
to $Ae_{i_{n-1}},$ but $Ae_{i_{n}} \rightarrow Ae_{i_{n-1}}
\rightarrow Ae_{i_{n}}$ is not equal to zero, hence $g$ can not be a
chain map.

c)The vertex $i_n$ is a degree one sink of $T_j$ but not $T_i$. Let
$i_{n+1}$ be the marked vertex of $T_i$ next to $i_n$. The
restrictions of $T_i$ and $T_j$ are arranged on the quiver of $A$ as
follows:

$$
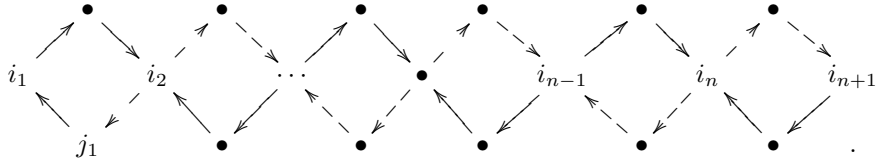
\xymatrix @=1pc{
&\bullet \ar[dr] &&\bullet\ar@{-->}[dr]&&\bullet\ar[dr]&&\bullet\ar@{-->}[dr]&&\bullet\ar[dr]& &\bullet\ar@{-->}[dr]& \\
i_1\ar[ur]&&i_2\ar@{-->}[ur]\ar@{-->}[dl]&&\cdots\ar[ur]\ar[dl]&&\bullet\ar@{-->}[ur]\ar@{-->}[dl]&&i_{n-1}\ar[dl]\ar[ur]&&i_n\ar@{-->}[dl] \ar@{-->}[ur] && i_{n+1}\ar[dl] \\
&j_1\ar[ul]&&\bullet\ar[ul]&&\bullet\ar@{-->}[ul]&&\bullet\ar[ul]&&\bullet\ar@{-->}[ul]& &\bullet\ar[ul]&. \\
}
$$

Let us construct a nonzero morphism from
$\emph{Hom}_{D^b(A)}(T_j,T_i[-1])$. This morphism is constructed the
same way as in the case when neither the restriction of $T_i$ with
respect to $T_j$ belongs to the restriction of $T_j$ with respect to
$T_i$ nor the restriction of $T_j$ with respect to $T_i$ belongs to
the restriction of $T_i$ with respect to $T_j$, and $i_1 \neq j_1,
i_n \neq j_n$. Let us start from the morphism $Ae_{i_{2}}
\rightarrow Ae_{i_{1}}$ and let us construct the morphism as
described earlier until we reach $Ae_{i_{n}} \rightarrow
Ae_{i_{n-1}}$, which was chosen such that the composition
$Ae_{i_{n-1}} \rightarrow Ae_{i_{n-2}} \rightarrow Ae_{i_{n-1}}$ is
equal to the composition $Ae_{i_{n-1}} \rightarrow Ae_{i_{n}}
\rightarrow Ae_{i_{n-1}}$ with the opposite sign, we can see that
the composition $Ae_{i_{n}} \rightarrow Ae_{i_{n-1}} \rightarrow
Ae_{i_{n}}$ is not equal to zero, choose $Ae_{i_{n}} \rightarrow
Ae_{i_{n+1}}$ such that $Ae_{i_{n}} \rightarrow Ae_{i_{n+1}}
\rightarrow Ae_{i_{n}}$ is equal to the composition $Ae_{i_{n}}
\rightarrow Ae_{i_{n-1}} \rightarrow Ae_{i_{n}}$ with the opposite
sign. Thus we have constructed a nonzero chain map from
$\emph{Hom}_{D^b(A)}(T_j,T_i[-1])$.

5) At last we can consider the case $i_1=j_1$, $i_n = j_n.$ The
following options are possible: a) $i_1$ is a degree one source of
$T_i,$ but not $T_j$, $i_n$ is a degree one sink of $T_i,$ but not
$T_j$; b) $i_1$ is a degree one source of $T_i,$ but not $T_j$,
$i_n$ is a degree one sink of both $T_i$ and $T_j$; c) $i_1$ is a
degree one source of both $T_i$ and $T_j$, $i_n$ is a degree one
sink of $T_i,$ but not $T_j$; d) $i_1$ is a degree one source of
$T_i,$ but not $T_j$, $i_n$ is a degree one sink of $T_j,$ but not
$T_i$.

a) The vertex $i_1$ is a degree one source of $T_i,$ but not $T_j$,
$i_n$ is a degree one sink of $T_i,$ but not $T_j$. Let $j_{0}$ be
the marked vertex of $T_j$ next to $i_1$, à $j_{n+1}$ be the marked
vertex of $T_j$ next to $i_n$. A nonzero chain map from the
restriction of $T_i$ to the restriction of $T_j$ shifted by $-1$,
should have a nonzero component from $Ae_{i_2}$ to $Ae_{i_1},$ but
$Ae_{i_{1}} \rightarrow Ae_{i_{2}} \rightarrow Ae_{i_{1}}$ is not
equal to zero, hence $\emph{Hom}_{D^b(A)}(T_i,T_j[-1])=0.$ A nonzero
chain map from the restriction of $T_j$ to the restriction of $T_i$
shifted by $-1$, should have a nonzero component from $Ae_{i_n}$ to
$Ae_{i_{n-1}},$ but $Ae_{i_{n}} \rightarrow Ae_{i_{n-1}} \rightarrow
Ae_{i_{n}}$ is not equal to zero, hence
$\emph{Hom}_{D^b(A)}(T_j,T_i[-1])=0.$

b) The vertex $i_1$ is a degree one source of $T_i,$ but not $T_j$,
$i_n$ is a degree one sink of both $T_i$ and $T_j$. Let $j_{0}$ be
the marked vertex of $T_j$ next to $i_1$. As in the previous case
$\emph{Hom}_{D^b(A)}(T_i,T_j[-1])=0.$ A nonzero chain map from the
restriction of $T_j$ to the restriction of $T_i$ shifted by $-1$,
should have a nonzero component from $Ae_{i_n}$ to $Ae_{i_{n-1}},$
but $Ae_{i_{n}} \rightarrow Ae_{i_{n-1}} \rightarrow Ae_{i_{n}}$ is
not equal to zero, hence $\emph{Hom}_{D^b(A)}(T_j,T_i[-1])=0.$

c) The vertex $i_1$ is a degree one source of both $T_i$ and $T_j$,
$i_n$ is a degree one sink of $T_i$ but not $T_j$. For reasons
similar to the case (a) we have
$\emph{Hom}_{D^b(A)}(T_j,T_i[-1])=0.$ A map from the restriction of
$T_i$ to the restriction of $T_j$ shifted by $-1$, should have a
nonzero component from $Ae_{i_2}$ to $Ae_{i_1},$ but $Ae_{i_{1}}
\rightarrow Ae_{i_{2}} \rightarrow Ae_{i_{1}}$ is not equal to zero,
hence $\emph{Hom}_{D^b(A)}(T_i,T_j[-1])=0.$

d) The vertex $i_1$ is a degree one source of $T_i,$ but not $T_j$,
$i_n$ is a degree one sink of $T_j,$ but not $T_i$. Let $j_{0}$ be
the marked vertex of $T_j$ next to $i_1$, and $i_{n+1}$ be the
marked vertex of $T_i$ next to $i_n$. Let us start the construction
of a nonzero morphism from $T_j$ to $T_i[1]$ with an arbitrary
nonzero morphism $Ae_{j_0}\rightarrow Ae_{i_1},$ choose a morphism
$Ae_{i_2}\rightarrow Ae_{i_1}$ such that $Ae_{i_1}\rightarrow
Ae_{j_0}\rightarrow Ae_{i_1}$ is equal to the composition
$Ae_{i_1}\rightarrow Ae_{i_2}\rightarrow Ae_{i_1}$ with the opposite
sign, as before let us construct a morphism from $T_j$ to $T_i[1]$
in such a manner that its composition with the differentials of
$T_j$ and $T_i[1]$ is equal to zero, let us finish the construction
choosing a morphism $Ae_{i_n}\rightarrow Ae_{i_{n+1}}$ such that the
composition $Ae_{i_n}\rightarrow Ae_{i_{n-1}}\rightarrow Ae_{i_n}$
is equal to the composition $Ae_{i_n}\rightarrow
Ae_{i_{n+1}}\rightarrow Ae_{i_n}$ with the opposite sign, thus
$\emph{Hom}_{D^b(A)}(T_j,T_i[-1]) \neq 0.$~\hfill\(\Box\)

\begin{lem}

Let the intersection of $T_i$ and $T_j$ consist of one connected
component, which contains more than one vertex. The diagrams of
$T_i$ and $T_j$ intersect in such a way that both of the end
vertices of the intersection are sinks, then
$\emph{Hom}_{D^b(A)}(T_i,T_j[-1])=0=\emph{Hom}_{D^b(A)}(T_j,T_i[-1])$
iff one of the following conditions hold: 1) neither the restriction
of $T_i$ with respect to $T_j$ belongs to the restriction of $T_j$
with respect to $T_i$ nor the restriction of $T_j$ with respect to
$T_i$ belongs to the restriction of $T_i$ with respect to $T_j$; 2)
the diagrams of $T_i$ and $T_j$ have a coinciding vertex of degree
one.

\end{lem}

\begin{lem}

Let the intersection of $T_i$ and $T_j$ consist of one connected
component, which contains more than one vertex. The diagrams of
$T_i$ and $T_j$ intersect in such a way that both of the end
vertices of the intersection are sources, then
$\emph{Hom}_{D^b(A)}(T_i,T_j[-1])=0=\emph{Hom}_{D^b(A)}(T_j,T_i[-1])$
iff one of the following conditions hold: 1) neither the restriction
of $T_i$ with respect to $T_j$ belongs to the restriction of $T_j$
with respect to $T_i$ nor the restriction of $T_j$ with respect to
$T_i$ belongs to the restriction of $T_i$ with respect to $T_j$; 2)
the diagrams of $T_i$ and $T_j$ have a coinciding vertex of degree
one.

\end{lem}

\textbf{Proof.} The proof of lemmas 6 and 7 is completely analogous
to that of lemma 5 and is left to the reader.~\hfill\(\Box\)

\begin{zam}

From lemmas 5, 6, 7 we can see that regardless of the type of the
intersection if two diagrams have a common sink or source of degree
one then the direct sum of the corresponding complexes is a
partially tilting complex.

\end{zam}
Let us now consider the case when at least one of the complexes is a
stalk complex of a projective module. Note that if $P$ is an
indecomposable stalk complex of a projective module concentrated in
degree 0 and $T_i$ is concentrated in 0 and 1 then the condition
$\emph{Hom}_{D^b(A)}(P,T_i[-1])=0$ holds automatically.

\begin{lem}

Let $T_i$ be the projective presentation of some module, such that
$T_i$ is an indecomposable two-term partial tilting complex and let
$P$ be an indecomposable stalk complex of a projective module
concentrated in degree 0, then $\emph{Hom}_{D^b(A)}(T_i,P[-1])=0$
iff one of the following conditions is satisfied: 1) the vertex
corresponding to $P$ coincides with a source of degree one of $T_i$;
2) the vertex corresponding to $P$ does not belong to $T_i$ and
there is no $A$-cycle which contains the vertex corresponding to $P$
and a sink of degree one of $T_i$ such that this is the only vertex
of $T_i$ belonging to this $A$-cycle.

\end{lem}

\textbf{Proof.} We have to consider the following cases: the vertex
corresponding to $P$ does not belong to $T_i,$ the vertex
corresponding to $P$ coincides with a source of $T_i,$ the vertex
corresponding to $P$ coincides with a sink of $T_i,$ the vertex
corresponding to $P$ coincides with an unmarked vertex of $T_i.$

1) The vertex corresponding to $P$ does not belong to $T_i.$ If the
vertex corresponding to $P$ does not belong to any $A$-cycle
containing vertices of $T_i$ then we clearly have
$\emph{Hom}_{D^b(A)}(T_i,P[-1])=0$; if the vertex corresponding to
$P$ and some substring of $T_i$ belong to the same $A$-cycle, then
the only morphism, from the sink of $T_i$ belonging to the same
$A$-cycle to $P$ does not induce a chain map, hence
$\emph{Hom}_{D^b(A)}(T_i,P[-1])=0$. If the vertex corresponding to
$P$ and just a source of $T_i$ belong to the same $A$-cycle, then
there is no nonzero maps from $(T_i)^1$ to $P$. If the vertex
corresponding to $P$ and just a sink of $T_i$ belong to the same
$A$-cycle, then the morphism from this sink of $T_i$ to $P$ induces
a nonzero chain map.

2) Let the vertex corresponding to $P$ coincide with the degree one
source of $T_i$. Let $i_2$ be the marked vertex of $T_i$ next to
$i_1$, then $\emph{Hom}_{D^b(A)}(T_i,P[-1])=0$ iff
$\emph{Hom}_{D^b(A)}(Ae_{i_1}\rightarrow Ae_{i_2},Ae_{i_1}[-1]=0$.
This is true since the only morphism $Ae_{i_2}\rightarrow Ae_{i_1}$
does not induce a chain map. Let the vertex corresponding to $P$
coincide with the degree two source $i_2$ of $T_i$. Let $i_1$, $i_3$
be the neighboring marked vertices of $T_i$, let us construct a
nonzero morphism $T_i\rightarrow P[-1]$ as follows. Choose an
arbitrary nonzero restriction of this map to $Ae_{i_1}\rightarrow
Ae_{i_2}$, the restriction to $Ae_{i_3}\rightarrow Ae_{i_2}$ is
chosen in such a manner that the composition $Ae_{i_2} \rightarrow
Ae_{i_1}\rightarrow Ae_{i_2}$ is equal to the composition $Ae_{i_2}
\rightarrow Ae_{i_3}\rightarrow Ae_{i_2}$ with the opposite sign,
set all other components to be zero. Thus,
$\emph{Hom}_{D^b(A)}(T_i,P[-1]) \neq 0$.

3) If the vertex corresponding to $P$ coincides with a sink $i_1$ of
$T_i$, then the map $Ae_{i_1}\rightarrow Ae_{i_1}$ which has the
socle of $Ae_{i_1}$ as its image induces a nonzero chain map and
$\emph{Hom}_{D^b(A)}(T_i,P[-1]) \neq 0$.

4) Let us consider the case when the vertex $k$ corresponding to $P$
coincides with an unmarked vertex of $T_i$. Let $i_1$, $i_2$ be the
source and the sink of the substring of $T_i$ containing $k.$
$\emph{Hom}_{D^b(A)}(T_i,P[-1])=0$ iff
$\emph{Hom}_{D^b(A)}(Ae_{i_1}\rightarrow Ae_{i_2},Ae_{k}[-1]=0$. It
is clear that a nonzero morphism $Ae_{i_2}\rightarrow Ae_{k}$
induces a chain map since the composition $Ae_{i_1}\rightarrow
Ae_{i_2}\rightarrow Ae_{k}$ is equal to zero.~\hfill\(\Box\)

\begin{lem}
Let $T_i$ be the projective presentation of some module, such that
$T_i$ is an indecomposable two-term partial tilting complex and let
$P$ be an indecomposable stalk complex of a projective module
concentrated in degree 1, then $\emph{Hom}_{D^b(A)}(T_i,P[1])=0$ iff
one of the following conditions is satisfied: 1) the vertex
corresponding to $P$ coincides with a sink of degree one of $T_i$;
2) the vertex corresponding to $P$ does not belong to $T_i$ and
there is no $A$-cycle which contains the vertex corresponding to $P$
and a source of degree one of $T_i$ such that this is the only
vertex of $T_i$ belonging to this $A$-cycle.

\end{lem}

\textbf{Proof.} The proof is analogous to that of lemma
8.~\hfill\(\Box\)

\begin{lem}

Let $P_i$, $P_j$ be two indecomposable stalk complex of projective
modules
$\emph{Hom}_{D^b(A)}(P_i,P_j[-1])=0=\emph{Hom}_{D^b(A)}(P_j,P_i[-1])$
iff one of the following conditions is satisfied: 1) the vertices
corresponding to $P_i$ and $P_j$ do not belong to the same $A$-cycle
2) the vertices corresponding to $P_i$ and $P_j$ belong to the same
$A$-cycle and $P_i$, $P_j$ are concentrated in the same degree.

\end{lem}

\section{Endomorphism rings}

In the previous section all two-term tilting complexes over a Brauer
tree algebra with multiplicity one were described as direct sums of
$n$ indecomposable partial tilting complexes which pairwisely
satisfy some conditions. In this section the endomorphism rings of
such complexes are described. It is well known that the endomorphism
rings are isomorphic to some Brauer tree algebras with multiplicity
one. Let $T$ be a tilting complex, to describe its endomorphism ring
it is sufficient to determine the partition of the vertices of
$\text{End}_{K^b(A)}(T)$ into $A$-cycles or equivalently which edges
of the Brauer tree of $\text{End}_{K^b(A)}(T)$ are incident to the
same vertex and the cyclic ordering of the edges incident to the
same vertex in the Brauer tree. The vertices of the quiver of
$\text{End}_{K^b(A)}(T)$ correspond to the indecomposable summands
of $T$, two vertices $i, j$ belong to the same $A$-cycle iff
$\text{Hom}_{K^b(A)}(T_i,T_j) \neq 0$ holds for the corresponding
summands $T_i, T_j$. This condition is easy to verify by the well
known formula due to Happel \cite{Ha2}: let $Q=(Q^r)_{r \in
\mathbb{Z}}, R=(R^s)_{s \in \mathbb{Z}} \in A\text{-}{\rm perf}$,
then
$$\sum_i (-1)^i {\rm dim}_K {\rm Hom}_{K^b(A)}(Q,R[i])=\sum_{r,s} (-1)^{r-s}{\rm dim}_K {\rm Hom}_{A}(Q^r,R^s).$$
Note that if ${\rm Hom}_{K^b(A)}(Q,R[i])=0, i \neq 0$ (for example,
in the case when $Q$ and $R$ are summands of a tilting complex) then
the left hand side of the formula becomes ${\rm dim}_K {\rm
Hom}_{K^b(A)}(Q,R).$

Recall that it is convenient to assume that each vertex belongs to
two $A$-cycles in the quiver of $A$, i.e. if some vertex has degree
2, then we assume that there is a formal loop, which is equal to the
other $A$-cycle passing through this vertex and this loop is the
second $A$-cycle containing this vertex.

\begin{predl}
Let $T_i,$ $T_j$ be projective presentations of nonisomorphic
indecomposable nonprojective $A$-modules such that their direct sum
is a partially tilting complex. $\text{Hom}_{K^b(A)}(T_i,T_j) \neq
0$ iff a degree one sinks of $T_i$ and $T_j$ belong to the same
$A$-cycle and these are the only vertices of $T_i$ and $T_j$,
belonging to this $A$-cycle, or a degree one sources of $T_i$ and
$T_j$ belong to the same $A$-cycle and these are the only vertices
of $T_i$ and $T_j$, belonging to this $A$-cycle.
\end{predl}

\textbf{Proof.} Let a degree one sinks of $T_i$ and $T_j$ belong to
the same $A$-cycle in such a manner that these are the only vertices
of $T_i$ and $T_j$, belonging to this $A$-cycle. This holds in the
following cases: the diagrams of $T_i$ and $T_j$ do not intersect, a
degree one sinks of $T_i$ and $T_j$ belong to the same $A$-cycle and
these are the only vertices of $T_i$ and $T_j$, belonging to this
$A$-cycle; $T_i$ and $T_j$ intersect and have a common sink of
degree one. The case of the sources is analogous.

1) Let $T_i$ and $T_j$ not have vertices belonging to the same
$A$-cycle, it is clear that ${\rm dim}_K {\rm
Hom}_{K^b(A)}(T_i,T_j)=0.$

2) Let $T_i$ and $T_j$ not intersect but have vertices belonging to
the same $A$-cycle $\Upsilon$. Let a substring $(i_1,...,i_2)$ of
$T_i$ and a substring $(j_1,...,j_2)$ of $T_j$ belong to $\Upsilon$,
the substrings do not intersect. It is clear that ${\rm dim}_K {\rm
Hom}_{K^b(A)}(T_i,T_j)={\rm dim}_K {\rm Hom}_{K^b(A)}(Ae_{i_2}
\rightarrow Ae_{i_1},Ae_{j_2} \rightarrow Ae_{j_1})=1-1-1+1=0.$ Now
let a substring $(i_1,...,i_2)$ of $T_i$ and an unmarked vertex of
$T_j$ belong to $\Upsilon$, we get ${\rm dim}_K {\rm
Hom}_{K^b(A)}(T_i,T_j)=0$, since the projective summands of the
components of $T_i$ and $T_j$ do not belong to the same $A$-cycle.
If $\Upsilon$ contains a degree one sink of $T_i$ and an unmarked
vertex of $T_j$; a degree one source of $T_i$ and an unmarked vertex
of $T_j$, then ${\rm dim}_K {\rm Hom}_{K^b(A)}(T_i,T_j)=0$ for the
same reason. Let a substring $(i_1,...,i_2)$ of $T_i$ and a degree
one sink $(j_1)$ of $T_j$ belong to $\Upsilon$. ${\rm dim}_K {\rm
Hom}_{K^b(A)}(T_i,T_j)={\rm dim}_K {\rm Hom}_{K^b(A)}(Ae_{i_2}
\rightarrow Ae_{i_1}, 0\rightarrow Ae_{j_1})=-1+1=0.$ The case when
a substring $(i_1,...,i_2)$ of $T_i$ and a degree one source $(j_1)$
of $T_j$ belong to $\Upsilon$ is similar to the previous one. Now
let a vertex $(i_1)$ of $T_i$ and a vertex $(j_1)$ of $T_j$ belong
to $\Upsilon$, both vertices are degree one sources. ${\rm dim}_K
{\rm Hom}_{K^b(A)}(T_i,T_j)={\rm dim}_K {\rm Hom}_{K^b(A)}(Ae_{i_1}
\rightarrow 0,Ae_{j_1} \rightarrow 0)=1.$ The case when $\Upsilon$
contains a degree one sink of both $T_i$ and $T_j$ is analogous.

3)Let the diagrams of $T_i$ and $T_j$ meet at one vertex $k.$ Let
$k$ be an unmarked vertex of $T_i$ and $T_j,$ then the marked
vertices of $T_i$ and $T_j$ belong to different $A$-cycles and there
are no nonzero morphisms between the corresponding projective
modules, ${\rm dim}_K {\rm Hom}_{K^b(A)}(T_i,T_j)=0.$ Let $k$ be a
degree one sink of both $T_i$ and $T_j,$ then ${\rm dim}_K {\rm
Hom}_{K^b(A)}(T_i,T_j)={\rm dim}_K {\rm
Hom}_{K^b(A)}(Ae_{i_2}\rightarrow Ae_k,Ae_{j_2}\rightarrow
Ae_k)=-1-1+2=0,$ where $i_2$ is the source of $T_i$ next to $k$,
$j_2$ is the source of $T_j$ next to $k$ and $i_2$ and $j_2$ belong
to different $A$-cycles. The case when $k$ is a degree one source of
both $T_i$ and $T_j$ is similar to the previous one.

Assume that the restriction of $T_i$ with respect to $T_j$ belongs
to the restriction of $T_j$ with respect to $T_i$ and the end points
of the intersection of $T_i$ and $T_j$ do not coincide with the
marked vertices of $T_j$. The corresponding subdiagrams are arranged
on the quiver of $A$ as follows:

4) Let $T_i$ and $T_j$ intersect in such a way that one of the end
points of the intersection is a sink and the other one is a source
(the intersection contains more than one vertex). Assume that the
restriction of $T_i$ with respect to $T_j$ belongs to the
restriction of $T_j$ with respect to $T_i$. The corresponding
subdiagrams are arranged on the quiver of $A$ as follows: (in the
case when $i_1 \neq j_1, i_n \neq j_n$):

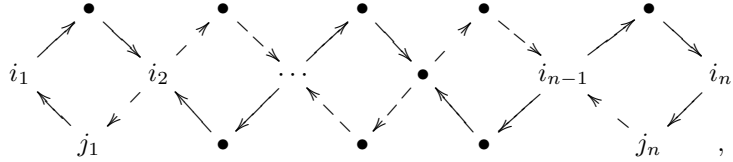
$$
\xymatrix @=1pc{
&\bullet \ar[dr] &&\bullet\ar@{-->}[dr]&&\bullet\ar[dr]&&\bullet\ar@{-->}[dr]&&\bullet\ar[dr]& \\
i_1\ar[ur]&&i_2\ar@{-->}[ur]\ar@{-->}[dl]&&\cdots\ar[ur]\ar[dl]&&\bullet\ar@{-->}[ur]\ar@{-->}[dl]&&i_{n-1}\ar[dl]\ar[ur]&&i_n\ar[dl] \\
&j_1\ar[ul]&&\bullet\ar[ul]&&\bullet\ar@{-->}[ul]&&\bullet\ar[ul]&&j_n\ar@{-->}[ul]&, \\
}
$$
where $i_1, i_2,...,i_n$ are the marked vertices of $T_i$ (which
belong to the restriction), the vertices with odd indices are
sources, the vertices with even indices are sinks; $j_1,
i_2,...,i_n-1,j_n$ are the marked vertices of $T_j$ (which belong to
the restriction), the vertices with odd indices are sources, the
vertices with even indices are sinks. Between the following modules
the Hom-space is one dimensional: $Ae_{i_k}$ and $Ae_{i_{k-1}},
Ae_{i_{k+1}}$ for $k=2,3,...,n-1$ and between $Ae_{i_1}$,
$Ae_{i_{2}},$ $Ae_{j_1}$, and also between $Ae_{i_{n-1}}$ and
$Ae_{i_{n-2}},$ and between $Ae_{i_{n}},$ $Ae_{i_{n-1}}$ and
$Ae_{j_{n}}$. As it was already mentioned, there are no nonzero
morphisms between the projective summands of the components of $T_i$
and $T_j$ which do not belong to the restrictions. It is easy to see
that ${\rm dim}_K {\rm Hom}_{A}(T_i^0,T_j^0)={\rm dim}_K {\rm
Hom}_{A}(T_i^0,T_j^1)={\rm dim}_K {\rm Hom}_{A}(T_i^1,T_j^0)={\rm
dim}_K {\rm Hom}_{A}(T_i^1,T_j^1)=n-1,$ hence using the formula we
get that ${\rm dim}_K {\rm Hom}_{K^b(A)}(T_i,T_j)=0$.

Now let $i_1=j_1$, $i_n \neq j_n,$ assume that $i_n$ belongs to the
restriction of $T_j$, the vertex $i_1$ is a degree one source of
both $T_i$ and $T_j$, then ${\rm dim}_K {\rm
Hom}_{A}(T_i^0,T_j^0)=n$, ${\rm dim}_K {\rm
Hom}_{A}(T_i^0,T_j^1)={\rm dim}_K {\rm Hom}_{A}(T_i^1,T_j^0)={\rm
dim}_K {\rm Hom}_{A}(T_i^1,T_j^1)=n-1,$ hence by the formula ${\rm
dim}_K {\rm Hom}_{K^b(A)}(T_i,T_j)=1$.

Let the vertex $i_1$ be a degree one source of $T_i$ but not $T_j$.
Let $j_0$ be the marked vertex of $T_j$ next to $j_1$. The
restrictions of $T_i$ and $T_j$ are arranged on the quiver of $A$ as
follows:

$$
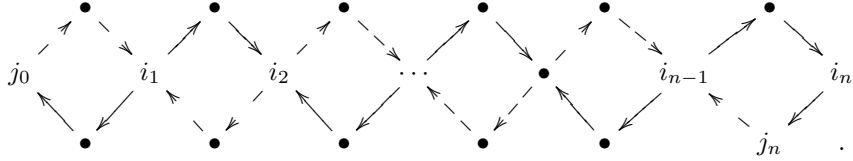
\xymatrix @=1pc{
&\bullet \ar@{-->}[dr] &&\bullet \ar[dr] &&\bullet\ar@{-->}[dr]&&\bullet\ar[dr]&&\bullet\ar@{-->}[dr]&&\bullet\ar[dr]& \\
j_0 \ar@{-->}[ur] &&i_1\ar[ur] \ar[dl]&&i_2\ar@{-->}[ur]\ar@{-->}[dl]&&\cdots\ar[ur]\ar[dl]&&\bullet\ar@{-->}[ur]\ar@{-->}[dl]&&i_{n-1}\ar[dl]\ar[ur]&&i_n\ar[dl] \\
&\bullet \ar[ul] &&\bullet\ar@{-->}[ul]&&\bullet\ar[ul]&&\bullet\ar@{-->}[ul]&&\bullet\ar[ul]&&j_n\ar@{-->}[ul]&. \\
}
$$

And we have ${\rm dim}_K {\rm Hom}_{A}(T_i^0,T_j^0)={\rm dim}_K {\rm
Hom}_{A}(T_i^0,T_j^1)=n$, ${\rm dim}_K {\rm
Hom}_{A}(T_i^1,T_j^0)={\rm dim}_K {\rm Hom}_{A}(T_i^1,T_j^1)=n-1,$
hence by the formula we get ${\rm dim}_K {\rm
Hom}_{K^b(A)}(T_i,T_j)=0$.

Now let $i_n=j_n$, $i_1 \neq j_1,$ assume that $i_1$ belongs to the
restriction of $T_j$, the vertex $i_n$ is a degree one sink of both
$T_i$ and $T_j$. Then ${\rm dim}_K {\rm Hom}_{A}(T_i^0,T_j^0)={\rm
dim}_K {\rm Hom}_{A}(T_i^0,T_j^1)={\rm dim}_K {\rm
Hom}_{A}(T_i^1,T_j^0)=n-1,$ ${\rm dim}_K {\rm
Hom}_{A}(T_i^1,T_j^1)=n$ hence by the formula ${\rm dim}_K {\rm
Hom}_{K^b(A)}(T_i,T_j)=1$.

Let the vertex $i_n$ be a degree one sink of $T_i$ but not $T_j$.
Let $j_{n+1}$ be the marked vertex of $T_j$ next to $i_n$. The
restrictions of $T_i$ and $T_j$ are arranged on the quiver of $A$ as
follows:

$$
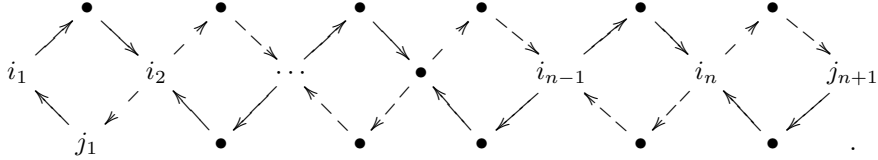
\xymatrix @=1pc{
&\bullet \ar[dr] &&\bullet\ar@{-->}[dr]&&\bullet\ar[dr]&&\bullet\ar@{-->}[dr]&&\bullet\ar[dr]& &\bullet\ar@{-->}[dr]& \\
i_1\ar[ur]&&i_2\ar@{-->}[ur]\ar@{-->}[dl]&&\cdots\ar[ur]\ar[dl]&&\bullet\ar@{-->}[ur]\ar@{-->}[dl]&&i_{n-1}\ar[dl]\ar[ur]&&i_n\ar@{-->}[dl] \ar@{-->}[ur] && j_{n+1}\ar[dl] \\
&j_1\ar[ul]&&\bullet\ar[ul]&&\bullet\ar@{-->}[ul]&&\bullet\ar[ul]&&\bullet\ar@{-->}[ul]& &\bullet\ar[ul]&. \\
}
$$

Then, ${\rm dim}_K {\rm Hom}_{A}(T_i^0,T_j^0)={\rm dim}_K {\rm
Hom}_{A}(T_i^0,T_j^1)=n-1,$ ${\rm dim}_K {\rm
Hom}_{A}(T_i^1,T_j^0)={\rm dim}_K {\rm Hom}_{A}(T_i^1,T_j^1)=n$
hence by the formula we get ${\rm dim}_K {\rm
Hom}_{K^b(A)}(T_i,T_j)=0$.

At last, let us consider the case $i_1=j_1$, $i_n = j_n.$ Let the
vertex $i_1$ be a degree one source of $T_i,$ but not $T_j$, $i_n$ a
degree one sink of $T_i,$ but not $T_j$. Let $j_{0}$ be the marked
vertex of $T_j$ next to $i_1$, and $j_{n+1}$ be the marked vertex of
$T_j$ next to $i_n$. Then ${\rm dim}_K {\rm
Hom}_{A}(T_i^0,T_j^0)={\rm dim}_K {\rm Hom}_{A}(T_i^0,T_j^1)={\rm
dim}_K {\rm Hom}_{A}(T_i^1,T_j^0)={\rm dim}_K {\rm
Hom}_{A}(T_i^1,T_j^1)=n$ hence by the formula we get ${\rm dim}_K
{\rm Hom}_{K^b(A)}(T_i,T_j)=0$.

Now let $i_1$ be a degree one source of $T_i,$ but not $T_j$, $i_n$
a degree one sink of both $T_i$ and $T_j$. Let $j_{0}$ be the marked
vertex of $T_j$ next to $i_1$. Then, ${\rm dim}_K {\rm
Hom}_{A}(T_i^0,T_j^0)={\rm dim}_K {\rm Hom}_{A}(T_i^0,T_j^1)={\rm
dim}_K {\rm Hom}_{A}(T_i^1,T_j^1)=n$, ${\rm dim}_K {\rm
Hom}_{A}(T_i^1,T_j^0)=n-1,$ hence by the formula we get ${\rm dim}_K
{\rm Hom}_{K^b(A)}(T_i,T_j)=1$.

Let the vertex $i_1$ be a degree one sink of both $T_i$ and $T_j$,
$i_n$ a degree one sink of $T_i$ but not $T_j$. Then ${\rm dim}_K
{\rm Hom}_{A}(T_i^0,T_j^0)={\rm dim}_K {\rm
Hom}_{A}(T_i^1,T_j^0)={\rm dim}_K {\rm Hom}_{A}(T_i^1,T_j^1)=n$,
${\rm dim}_K {\rm Hom}_{A}(T_i^0,T_j^1)=n-1,$ hence by the formula
we get ${\rm dim}_K {\rm Hom}_{K^b(A)}(T_i,T_j)=1$.

As it was mentioned before the following types of intersections can
occur: one of the end points of the intersection is a sink and
another is a source (this case we have just analysed) , both of the
end points of the intersection are sinks, both of the end points of
the intersection are sources. The cases when both of the end points
of the intersection are sinks or both of the end points of the
intersection are sources can be analysed similar to the case 4 and
is left to the reader.~\hfill\(\Box\)

\begin{predl}
Let $T_i$ be the projective presentation of an indecomposable
nonprojective $A$-modules, $P$ a stalk complex of an indecomposable
projective module concentrated in degree 0 and let $T_i \bigoplus P$
be a partial tilting complex. $\text{Hom}_{K^b(A)}(T_i,P) \neq 0$
iff a degree one source of $T_i$ and the vertex corresponding to $P$
belong to the same $A$-cycle and this source is the only vertex of
$T_i$ belonging to this $A$-cycle.
\end{predl}

\textbf{Proof.} Let $j$ denote the vertex corresponding to $P$. It
is clear that if $j$ and vertices of $T_i$ do not belong to the same
$A$-cycle, then $\text{Hom}_{K^b(A)}(T_i,P) = 0$. Let $j$ and some
vertices of $T_i$ belong to the same $A$-cycle $\Upsilon$, $j$ does
not belong to $T_i$, there are the following possibilities: 1) there
is a substring $(i_1,...,i_2)$ of $T_i$ which belongs to $\Upsilon$
2) there is just a degree one source $i$ of $T_i$ which belongs to
$\Upsilon$. In the first case ${\rm dim}_K {\rm
Hom}_{K^b(A)}(T_i,P)={\rm dim}_K {\rm
Hom}_{K^b(A)}(Ae_{i_2}\rightarrow Ae_{i_1},Ae_{j})=0$. In the second
case ${\rm dim}_K {\rm Hom}_{K^b(A)}(T_i,P)={\rm dim}_K {\rm
Hom}_{K^b(A)}(Ae_{i},Ae_{j})=1$.

Let now $j$ belong to the diagram of $T_i$, i.e. $j$ coincides with
a degree one source of $T_i.$ Denote by $i$ the marked sink of $T_i$
next to $j.$ Then ${\rm dim}_K {\rm Hom}_{K^b(A)}(T_i,P)={\rm dim}_K
{\rm Hom}_{K^b(A)}(Ae_{j}\rightarrow
Ae_{i},Ae_{j})=2-1=1$.~\hfill\(\Box\)

\begin{predl}
Let $T_i$ be the projective presentation of an indecomposable
nonprojective $A$-modules, $P$ a stalk complex of an indecomposable
projective module concentrated in degree 1 and let $T_i \bigoplus P$
be a partial tilting complex. $\text{Hom}_{K^b(A)}(T_i,P) \neq 0$
iff a degree one sink of $T_i$ and the vertex corresponding to $P$
belong to the same $A$-cycle and this sink is the only vertex of
$T_i$ belonging to this $A$-cycle.
\end{predl}

\textbf{Proof.} The proof is similar to the proof of proposition
2.~\hfill\(\Box\)

\begin{predl}
Let $P_i$, $P_j$ be stalk complexes of indecomposable projective
modules and let $P_i \bigoplus P_j$ be a partial tilting complex.
$\text{Hom}_{K^b(A)}(P_i,P_j) \neq 0$ iff $P_i$ and $P_j$ are
concentrated in the same degree and the vertices corresponding to
$P_i$ and $P_j$ belong to the same $A$-cycle.
\end{predl}

From the description of the Cartan matrix we see that there are two
types of $A$-cycles in $\text{End}_{K^b(A)}(T)$: $A$-cycles
corresponding to sources, to such an $A$-cycle belong all the
indecomposable partial tilting complexes with diagrams having a
degree one source on some fixed $A$-cycle $\Upsilon$ of algebra $A$
such that no other vertices of these diagrams belong to $\Upsilon$
and all the indecomposable stalk complexes of projective modules
concentrated in degree 0 such that the corresponding vertices belong
to $\Upsilon$; $A$-cycles corresponding to sinks, to such an
$A$-cycle belong all the indecomposable partial tilting complexes
with diagrams having a degree one sink on some fixed $A$-cycle
$\Upsilon$ of algebra $A$ such that no other vertices of these
diagrams belong to $\Upsilon$ and all the indecomposable stalk
complexes of projective modules concentrated in degree 1 such that
the corresponding vertices belong to $\Upsilon$.

Let us describe the cyclic ordering of the indecomposable summands
of $T$ belonging to the $A$-cycle corresponding to sources. With
this aim in view let us introduce a sort of cyclic lexicographic
order. Fix a vertex on $\Upsilon$ and set it to be the greatest (the
vertex with an arrow coming from the fixed vertex set to be less
than the fixed one, the next vertex even less and so on), order the
other vertices linearly. To a diagram with a degree one source on
$\Upsilon$ let us associate an ordered set of vertices as follows:
the first vertex is the degree one source on $\Upsilon$, after that
take all the marked vertices according to their order on the
diagram. To an indecomposable stalk complex of a projective module
associate a set consisting of the vertex which corresponds to that
module. Let us consider a usual lexicographic order on the sets of
vertices (except for we set the empty spot on an even position to be
the least and on an odd position the greatest): if the first vertex
of the set corresponding to $T_i$ is less than the first vertex of
the set corresponding to $T_j$, then $T_i<T_j$, if these vertices
coincide, then the second vertices of the sets corresponding to
$T_i$ and $T_j$ belong to the same $A$-cycle, the first vertex of
the sets belong to the same $A$-cycle, set the first vertex of the
sets to be the greatest (among the vertices) on this $A$-cycle, then
we can consider a linear order on this $A$-cycle as earlier, if the
second vertex of the set corresponding to $T_i$ is less than the
second vertex of the set corresponding to $T_j$, then $T_i<T_j$. If
all the vertices from the first to the $i$-th corresponding to the
sets of $T_i$ and $T_j$ coincide, and the $i+1$-st vertices of the
sets differ, then the $i$-th vertex and the $i+1$-st vertices belong
to the same $A$-cycle, setting the $i$-th vertex to be the greatest
among the vertices (the empty spot can be greater)we can compare the
$i+1$-st vertices. The cyclic ordering is glued from the linear one.

Recall that we identify the indecomposable summands of $T$ and the
edges in the Brauer tree of $\text{End}_{K^b(A)}(T)$.

\begin{predl}
In the Brauer tree of $\text{End}_{K^b(A)}(T)$ the cyclic ordering
of the edges incident to a vertex corresponding to an $A$-cycle of
sources coincides with the order introduced above.
\end{predl}

\textbf{Proof.} Let $\Upsilon$ be some $A$-cycle of algebra $A,$
denote the $A$-cycle of $\text{End}_{K^b(A)}(T)$ corresponding to
the sources belonging to $\Upsilon$ by $\Psi$. Assume that $r$
vertices belong to $\Psi$, then to determine the cyclic ordering of
the vertices belonging to $\Psi$ it is sufficient to construct $r$
morphisms between the corresponding summands of $T$ such that their
consecutive composition is not homotopic to zero. First let us
construct morphisms from $T_j$ to $T_i$, where $T_i<T_j$ then from
the least summand to the greatest. Let us construct $\alpha_{2,1}$,
assume that in the sets of the vertices corresponding to $T_2,$
$T_1$ the first $i$ vertices coincide and the $i+1$-st vertices are
different, $T_2 > T_1.$ Denote by $P_1,...,P_i$ the projective
modules corresponding to the first $i$ vertices of the sets of
$T_1,$ $T_2$, denote by $P_{i+1}$ the projective module
corresponding to the $i+1$-st vertex of $T_1,$ by $P_{i+2}$ the
projective module corresponding to the $i+1$-st vertex of $T_2$.

Let us set the morphism $\alpha_{2,1}:T_2 \rightarrow T_1$ on the
projective summands of the components of $T_1$ and $T_2$. On the
coinciding modules $P_1,...,P_i$ it is the identity morphism and
$\alpha_{2,1}|_{P_{i+2} \rightarrow P_{i+1}}$ is the multiplication
by the unique path in the quiver of $A$ between the corresponding
vertices. All other components are zero. Let us check that the map
obtained is a chain map. It is clear that the commutativity of the
square $$
\xymatrix { T_2^0 \ar[r] \ar[d]^{\alpha_{1,2}^0} & T_2^1 \ar[d]^{\alpha_{1,2}^1}  \\
T_1^0 \ar[r]  & T_2^1}
$$ follows from the commutativity of the squares consisting of
direct summands of $T_2^0, T_2^1, T_1^0, T_2^1$. The squares with
only identity or only zero morphisms are commutative. Let $i$ be
odd, the square $$
\xymatrix { P_i \ar[r] \ar[d]^{1} & P_{i+2} \ar[d]^{\alpha_{1,2}^1}  \\
P_i \ar[r]  & P_{i+1}}
$$ is commutative, since the path from the vertex corresponding to
$P_i$ to the vertex corresponding to $P_{i+1}$ passes the vertex
corresponding to $P_{i+2}$. If $P_{i+1}=0$, the square remains
commutative.

Let $i$ be even, the square $$
\xymatrix { P_{i+2}  \ar[r] \ar[d]^{\alpha_{1,2}^0} & P_i \ar[d]^{1}  \\
P_{i+1} \ar[r]  & P_i }
$$ is commutative, since the path from the vertex corresponding to
$P_{i+2}$ to the vertex corresponding to $P_{i}$ passes the vertex
corresponding to $P_{i+1}$. If $P_{i+2}=0$, the square remains
commutative. The second square containing $\alpha_{1,2}|_{P_{i+2}
\rightarrow P_{i+1}}$ is commutative since the modules it contains
belong to different $A$-cycles and the compositions are equal to
zero.

Let as construct the morphism from the least summand of $T$ to the
greatest $\alpha_{min, max}: T_{min} \rightarrow T_{max}$. The
projective modules corresponding to the first elements of the sets
of $T_{min}, T_{max}$ are denoted by $P_{min}, P_{max}$
respectively. If $P_{min} \neq P_{max}$ then the component
$\alpha_{min, max}| P_{min} \rightarrow P_{max}$ is set to be the
multiplication by the unique nonzero path, all other components are
set to be zero. If $P_{min} = P_{max}$, then the component
$\alpha_{min, max}| P_{min} \rightarrow P_{max}$ is set to be the
multiplication by the long nonzero path, i.e. the morphism with the
socle of $P_{max}$ as its image, all other components are set to be
zero. The commutativity of
$$
\xymatrix { T_{min}^0 \ar[r] \ar[d]^{\alpha_{1,2}^0} & T_{min}^1 \ar[d]^{\alpha_{1,2}^1}  \\
T_{max}^0 \ar[r]  & T_{max}^1} $$ follows from the commutativity of
$$
\xymatrix { P_{min} \ar[r] \ar[d]^{\alpha_{1,2}^0} & P_1 \ar[d]^{0}  \\
P_{max} \ar[r]  & P_2.} $$ Where $P_1, P_2$ are the modules
corresponding to the second elements of the sets of $T_{min},
T_{max}$ respectively (note that $P_1, P_2$ can be zero). If
$P_{min} \neq P_{max}$, then $P_{min}$ and $P_2$ belong to different
$A$-cycles, hence the composition $P_{min} \rightarrow P_{max}
\rightarrow P_2$ is equal to 0. If $P_{min} = P_{max}$, then the
composition $P_{min} \rightarrow P_{max} \rightarrow P_2$ is equal
to 0, since the image of  $P_{min} \rightarrow P_{max}$ is the socle
of $P_{max}$.

$A$-cycle $\Psi$ consists of $r$ vertices, let us check that the
consecutive composition of $r$ constructed morphisms is not
homotopic to 0. That is for any indecomposable summand of $T$ (say
$T_i$) the composition $\alpha_{i,i}: T_i \rightarrow T_i$ is not
homotopic to 0. Let $P_i$ be the module corresponding to the first
vertex of the set $T_i$, it's easy to see that $\alpha_{i,i}|_{P_i
\rightarrow P_i}$ is the multiplication by the long nonzero path,
i.e. the morphism with the socle of $P_{i}$ as its image, all other
components of $\alpha_{i,i}$ are zero. This morphism is not
homotopic to 0.~\hfill\(\Box\)

Let us describe the cyclic ordering of the indecomposable summands
of $T$ belonging to the $A$-cycle corresponding to sinks. It differs
from the ordering of the indecomposable summands of $T$ belonging to
the $A$-cycle of sources only by the fact that the sets
corresponding to the diagrams start with a degree one sink on
$\Upsilon$ and that the empty spot on the odd position is set to be
least and on the even position the greatest. Namely fix a vertex on
$\Upsilon$ and set it to be the greatest, order the other vertices
linearly. To a diagram with a degree one sink on $\Upsilon$ let us
associate an ordered set of vertices as follows: the first vertex is
the degree one sink on $\Upsilon$, after that take all the marked
vertices according to their order on the diagram. To an
indecomposable stalk complex of a projective module associate a set
consisting of the vertex which corresponds to that module. Let us
consider a lexicographic order on the sets of vertices (except for
we set the empty spot on an odd position to be the least and on an
even position the greatest) as before, glue the cyclic ordering from
this linear one. The proof of the next proposition is similar to the
proof of proposition 5.

\begin{predl}

In the Brauer tree of $\text{End}_{K^b(A)}(T)$ the cyclic ordering
of the edges incident to a vertex corresponding to an $A$-cycle of
sinks coincides with the order introduced above.

\end{predl}

\begin{predl}

Over any Brauer tree algebra $A$ there exists a two-term tilting
complex $T$ such that the algebra $\text{End}_{K^b(A)}(T)$ is
isomorphic to the Brauer star algebra.

\end{predl}

\textbf{Proof.} Fix some $A$-cycle $\Upsilon$ of algebra $A$. To any
vertex $x$ of algebra $A$ let us associate an indecomposable
two-term partial tilting complex $T_x$ such that the sum over all
the vertices is the desired $T$. All the summands of $T$ will belong
to an $A$-cycle of $\text{End}_{K^b(A)}(T)$, corresponding to
sources. If the vertex $x$ belong to $\Upsilon$, then $T_x$ is the
stalk complex of the projective module corresponding to $x$,
concentrated in 0. If the vertex $x$ does not belong to $\Upsilon$,
consider a diagram such that one of its end points is $x,$ and the
other end point is some vertex $y$ belonging to $\Upsilon$ such that
$y$ is the only vertex of the diagram belonging to $\Upsilon$ and
$y$ is a source. It is clear that since the Brauer graph of $A$ is a
tree, this diagram exists and is unique, $T_x$ is the indecomposable
complex corresponding to this diagram. From the construction we get
that $T:=\bigoplus_{x \in A}T_x$ is a tilting complex: firstly, $T$
contains the required number of nonisomorphic direct summands,
secondly, $T$ is partially tilting. If the diagrams corresponding to
different vertices $x$ and $y$ do not intersect, use lemma 2, if
they intersect then remark 2, if one of the vertices belong to
$\Upsilon$, use lemma 8, if both belong to $\Upsilon$, use lemma
10.~\hfill\(\Box\)

\begin{zam}

It is clear that similarly to the construction from proposition 7 we
could construct a tilting complex $T$ such that
$\text{End}_{K^b(A)}(T)$ is isomorphic to the Brauer star algebra
and all summands of $T$ belong to an $A$-cycle
$\text{End}_{K^b(A)}(T)$ corresponding to sinks.

\end{zam}

\begin{zam}
Over a Brauer tree algebra associated to a Brauer tree with $n$
edges there are exactly $2(n+1)$ nonisomorphic basic two-term
tilting complexes $T$ such that $\text{End}_{K^b(A)}(T)$ is
isomorphic to a Brauer star algebra. Each of $n+1$ cycles of $A$ can
generate a complex $T$ such that all its summands belong to
$A$-cycle of $\text{End}_{K^b(A)}(T)$ corresponding to sinks or
sources.

\end{zam}

\end{document}